\documentclass[12pt]{article}
\usepackage{graphicx}

\font\tendb=msbm10 at 12pt
\font\sevendb=msbm10 at 9pt
\font\fivedb=msbm10 at 7pt
\newfam\dbfam
\textfont\dbfam=\tendb
\scriptfont\dbfam=\sevendb
\scriptscriptfont\dbfam=\fivedb
\def\db{\fam\dbfam\tendb}
\font\eufm=eufm10\font\eufms=eufm10\font\eufmss=eufm10\newfam\eufam
\textfont\eufam=\eufm\scriptfont\eufam=\eufms\scriptscriptfont\eufam=\eufmss

\font\tendbb=msbm10 at 12pt
\font\sevendbb=msbm7 at 9pt
\font\fivedbb=msbm5 at 6pt
\newfam\dbbfam
\textfont\dbbfam=\tendbb
\scriptfont\dbbfam=\sevendbb
\scriptscriptfont\dbbfam=\fivedbb

 \def \Z {{\db Z}}
 \def \R {\hbox{\db R}}
 
 \def \C {\hbox{\db C}}
 \def \N {{\db N}}
 \def \S {S^{3}}

\font\tenMmm=eusm10 at 12pt
\newfam\Mmmfam
\textfont\Mmmfam=\tenMmm

\def\illu #1 by #2 (#3){
  \vbox to #2{
    \hrule width #1 height 0pt depth 0pt
    \vfill
    \special{illustration #3} 
    }
  }

\begin{document}

\null
\vspace{4cm}

\begin{center}  {\large {\bf  Quantum invariants and finite group actions on three-manifolds }}\\
 Nafaa Chbili\footnote{This work was partially completed during my stay at  the Research Institute for Mathematical Sciences, Kyoto University, Japan.
  I would like to express my thanks and gratitude to the Mitsubishi Foundation and the Sumitomo Foundation for their financial support.
  I am  also grateful to Hitoshi Murakami for his kind hospitality.}
\end{center}
\vspace{10mm}

 \begin{footnotesize}
               {\bf Abstract.} A 3-manifold $M$ is said to be $p$-periodic ($p\geq 2$ an integer) if and only if the finite
                cyclic group of order $p$ acts on $M$ with a circle as the set of fixed points. This paper provides
                 a criterion for periodicity of rational homology three-spheres. Namely, we give a necessary condition
                 for a rational homology three-sphere to be periodic with a prime period. This condition is given in
                 terms of the quantum $SU(3)$ invariant. We also discuss similar results for the Murakami-Ohtsuki-Okada
                 invariant.

               \end{footnotesize}

\begin{center}{\sc Introduction}\end{center}
The last decade has seen many new invariants introduced to low dimensional topology. In the case
of three-manifolds, these invariants were first proposed  by E. Witten \cite{Wi} who used the Feynman path integral based on the Chern-Simons gauge theory in order  to prove the existance of an invariant of oriented
3-manifolds associated with  each compact Lie group. The first rigorous construction of such invariants (called quantum invariants) was by   Reshetikhin and Turaev \cite{RT}  in the case of the Lie group $SU(2)$. Various   approaches to the construction of the quantum $SU(2)$ invariant  followed. Particularly, Lickorish \cite{Li2}  gave an elementary construction of this invariant using
the skein theory associated to the Kauffman bracket  \cite{Ka}  evaluated at primitive roots of unity of order
$4r$. Blanchet, Habegger, Masbaum and Vogel \cite{BHMV} then, extended Lickorish's approach to define the invariant
at roots of unity of  order $2p$. On the lines  of Lickorish's work, Ohtsuki and Yamada \cite{OY} defined skein linear
theory for $SU(3)$ using Kuperberg's skein relations  \cite{Ku} for trivalent oriented graphs. Via  this skein
theory, they constructed magic  elements to  define the quantum $SU(3)$ invariant of 3-manifolds at primitive roots of unity of order $6r$. Their work was extended by Miyazawa and Okamoto  \cite{MO} applying it  to primitive roots of unity of order $3r$. A more general construction was proposed  by Yokota who defined a  linear skein theory for quantum
$SU(N)$ invariants \cite{Yo}.\\
It is worth mentioning that  we have no idea about the power and the
preciseness of the above mentioned invariants  in detecting the topology of
 a given 3-manifold, despite the large number of approaches proposed to construct these invariants. In some sense the topological meaning of the quantum   invariants remains mysterious and still  far from being completely  understood.
The main purpose of this  paper is  to investigate the behaviour of the quantum $SU(3)$ invariant in the case where the manifold presents some geometric
properties (symmetries), discussing  whether or not
 the symmetry of the given manifold is reflected on its invariants.
 Ultimately,  by this we hope to   shed some light on the behaviour of
 the quantum invariants as well as   to  provide a criterion for periodicity of 3-manifolds.\\
Let $M$ be a 3-manifold (all the manifolds considered in this paper
are oriented  compact connected and  without boundary) and let $G$ be
the finite cyclic group $\Z/p\Z$. One of the most intriguing questions in
low dimensional topology is whether the group $G$ acts on $M$. Some
 progress in solving  this problem has been made  using classical
  algebraic and geometric tools. Nevertheless,  many related questions
   are still open. In this paper we shall  focus  our interests on  the
    case where $G$ acts on $M$ with a circle as the set of fixed points.
    A manifold with such an action is said to be $p$-periodic.
    P. Gilmer \cite{Gi}, studied the Witten-Reshitikhin-Turaev $SU(2)$
    invariant of periodic 3-manifolds as well as the case where the
    action
    of $G$ has no fixed points. On the other hand,   we used  different
    techniques in  \cite{Ch1}  to give a
necessary condition for a
rational homology 3-sphere to be periodic, this condition is given in
 terms of  the quantum $SU(2)$ invariant. Recently, other similar results
  were obtained by Bartoszynska, Gilmer and Przytycki \cite{PB},
  Chen and Le \cite{CL}.\\
  The present paper introduces a new criterion for periodicity of 3-manifolds using the quantum $SU(3)$ invariant. For technical reasons we will only deal  with the case of  rational homology 3-spheres. We also
prove similar results for the Murakami-Ohtsuki-Okada invariant (MOO for
short) \cite{MOO}. This invariant is known as  the simplest quantum
 invariant.\\
 This paper is organized as follows, in paragraph 1 we
  present our results. In paragraph 2 we explain the relationship between
  periodic links and periodic 3-manifolds. A review of  linear skein theory and definition of the quantum  $SU(3)$ invariant  are  given  in paragraph 3. The  proof of our
main theorem will be given  in paragraph 4.
 Further speculations are discussed in the last two paragraphs.\\
\begin{center}{\sc 1-Results and Application}\end{center}

Let $r \geq 4$ be  an integer and $M$ a 3-manifold. For
 $k \in \N ^{\star}$, let  $\Phi_{k}(A)$ denote  the cyclotomic polynomial
  of order $k$, $A$ here is an  indeterminate, and  $[k]$  denote the
  quantum integer $\frac{A^{3k}-A^{-3k}}{A^{3}-A^{-3}}$. Here we consider
  the quantum $SU(3)$ invariant defined at roots of unity of order $ 3r$.
  This invariant will be  denoted by ${\cal I}_{r}$ (see the definition in
  paragraph 3). Let $\Lambda_{r}= \Z[A^{\pm 1}]_{/\Phi_{3r}(A)}$, we shall
  explain later that the invariant
${\cal I}_{r}$ is an element of $\Lambda_{r}[\displaystyle\frac{1}{3r}]$. Consequently, if
$p$ is coprime to $3r$, it is possible
 to consider ${\cal I}_{r}$ with coefficients reduced modulo $p$.
  Throughout the rest of this paper, if $M$ is a periodic manifold then
   $\overline M$ denotes the quotient space. Recall
   here that if $M$ is  a rational homology sphere then $\overline M$
   is also a rational homology sphere. For $r\geq 4$ we will  denote by  $G_{r}$ the element of $\Lambda_{r}[\displaystyle\frac{1}{3r}]$ given by the following formula:
$$G_{r}(A)=A^{-36}\frac{ (\displaystyle\sum_{r=0}^{r-1}A^{6k^{2}})^{2}(\displaystyle\sum_{r=0}^{3r-1}A^{2k^{2}})^{2}}{3r^{2}}.$$

{\bf Theorem 1.1} {\sl Let $p$ be  a prime and $M$ a rational homology 3-sphere. If $M$ is $p$-periodic then for
all odd
 $r\geq 4$ such that $p$ and $3r$ are coprime,  the following congruence holds:
$$ {\cal I }_{r}(M)\equiv ({\cal I }_{r}(\overline M ))^{p} (G_{r}(A))^{\alpha} , \mbox{ modulo  }p ,[3]^{p}-[3] ,$$
where $\alpha$ is an integer.}\\

The congruence given by this theorem holds in the ring
 $\Lambda_{r}[\displaystyle\frac{1}{3r}]$. In this ring we have
 no idea about  how large the ideal generated by $p$ and  $ [3]^{p}-[3]$ is. However, if   one considers the case $p$ congrunent to $\pm 1$ modulo $r$, then we can show that working modulo $p$ and  $ [3]^{p}-[3]$ is equivalent to working modulo $p$. This is explained by the following corollary written  in the case of cyclic branched coverings of the 3-sphere.\\

{\bf Corollary 1.2} {\sl Let $p\neq 3$ be a  prime and $M$ the $p$-fold cyclic covering of $\S$  branched along a knot $K$. For all odd  $r\geq 4$ such that   $p \equiv \pm 1 \mbox { modulo  }r$,  the following congruence holds:
$$ {\cal I }_{r}(M)\equiv
 (G_{r}(A))^{\alpha}  \mbox{ mod }p  ,$$
where $\alpha$ is an integer.}\\

{\bf Remark 1.3} It should be  mentioned here that the condition of
 theorem 1.1 gives no obstruction to the  periodicity of a given rational
 homology 3-sphere. However, if we restrict the condition to cyclic
 branched coverings of $\S$ as in corollary 2.2, the condition takes
 a nice form and turns out to be easy to handle. Indeed, this condition
 is not trivial, in order  to prove this we can examine the case $r=5$. In this case $G_{5}(A)=A^{-36}$, hence modulo $p$,  ${\cal I }_{5}$ is a power of $A$. Thus,
given a rational homology 3-sphere,  we  just need
 to verify if $ {\cal I }_{r}(M)$ is  not a power of $A$ in order
 to conclude that $M$ is not the $p$-fold branched covering of $\S$(for appropriate values of $p$). This is  illustrated by the following example. Let
$L(2,1)$ be the lens space of type $(2,1)$. Put $r=5$, the cyclotomic polynomial of order 15 is given by $\Phi_{15}(A)=
1-A+A^{3}-A^{4}+A^{5}-A^{7}+A^{8}$. According to the general formula for the $SU(3)$ invariants of lens spaces obtained in \cite{MO}, we have the following:
$${\cal I}_{5}(L(2,1))=1-A-A^{2}+A^{3}-A^{4}+A^{5}-A^{7}.$$
The  powers of $A$ in the ring $\Lambda_{5}[\displaystyle\frac{1}{15}]$ are given by  the list below:
$$
\begin{array}{rl}
&1,\;A,\; A^{2},\; A^{3},\; A^{4},\; A^{5},\; A^{6}, \;A^{7},\\
&-1+A-A^{3}+A^{4}-A^{5}+A^{7},\\
&A^2-A^6-1-A^3+A^7,\\
&-1-A^5,\;-A-A^6,\;-A^2-A^7,\\
&1-A-A^{4}+A^{5}-A^{7},\\
&1-A^{2}+A^{3}-A^{4}+A^{6}-A^{7}.
\end{array}
$$
Obvisouly, ${\cal I}_{5}(L(2,1))$ does not belong to the list above. Consequently, the lens space  $L(2,1)$ is not the $p$-fold cyclic covering of the three-sphere branched along a knot $K$ for any prime $p$ congruent to $\pm 1$ modulo 5.


{\bf Proof of Corollaray 1.2.} If $M$ is the $p$-fold cyclic covering of $\S$  branched along a knot $K$, then $M$ is $p$-periodic and its  quotient is  $\S$.   Recall that the invariant ${\cal I }_{r}$ is normalized in such a way that   ${\cal I }_{r}(\S)=1$, for all $r \geq 4$.
 If $p\neq 3$ and congruent to  $\pm 1$ modulo $r$ then $p$ and $3r$ are coprime. Thus, we have modulo $p$:

$$\begin{array}{rl}
 [3]^{p}-[3]\equiv& (1+A^{6}+A^{-6})^{p}-(1+A^{6}+A^{-6})\\
\equiv &A^{6p}+A^{-6p}-A^{6}-A^{-6}
\equiv 0.
\end{array}
$$
The last equality holds because if we are working modulo  $\Phi_{3r}(A)$ then $A^{3r}=1$. Hence if $p\equiv 1$ modulo $r$ then $3r$ divides $6(p-1)$. Therefore  $A^{6(p-1)}=1$, that is
$A^{6p}=A^{6}$ and $A^{-6p}=A^{-6}$, same argument for  $p\equiv -1$ modulo $r$.\\

\begin{center}{\sc 2--Surgery presentation of periodic 3-manifolds }\end{center}

Let $M$ be a 3-manifold. It is well known from Lickorish's
 work \cite{Li1} that $M$ may be obtained from the 3-sphere $\S$ by
 surgery along a framed link $L$ in $\S$. Such a link is called a surgery
  presentation of $M$. Kirby \cite{Ki} introduced elementary moves
   which suffice for moving from one surgery presentation of $M$ to another.
    This section deals  with surgery presentation of periodic manifolds. We here recall briefly some defintions  and results, (see \cite{PS} for  details). \\

{\bf Definition 2.1} {\sl {
Let $p\geq 2$ an integer, a link $L$ of $\S$ is said to be $p$-periodic if and only if there exists an orientation preserving  auto-diffeomorphism  of  $\S$
such that:\\
1- Fix($h$) is homeomorphic to the  circle $S^{1}$,\\
2- the link  L is disjoint from  Fix($h$),\\
3- $h$ is of order $p$,\\
4- $h(L)=L$.\\
If $L$ is  periodic we will denote the quotient link by $\overline L$.\\}}

If we  consider $\S= \{(z_{1},z_{2}) \in \C ^{2}; |z_{1}|^{2}+|z_{2}|^{2}=1\}$,
then by the positive solution of the Smith conjecture
\cite{BM} we know that $h$ is topologically conjugate to the standard rotation given by:\\
 $$
\begin{array}{cccl}
\varphi_{p}:& S^{3} & \longrightarrow & S^{3} \\
   & (z_{1},z_{2}) & \longmapsto & (e^{\frac{2i\pi}{p}}z_{1},z_{2}).
\end{array} $$
The set of fixed points for  this diffeomorphism is the circle $\Delta = \{(0,z_{2}) \in \C ^{2}; |z_{2}|^{2}=1\}$. If $L$ is a  $p$-periodic link we denote by $lk(L,\Delta)$ the linking number of $L$ and the circle $\Delta$. If we identify $\S$ with $\R ^{3}\cup {\infty}$, $\Delta$ may be seen as  as the
the standard $z$-axis. For example,  the trefoil knot $3_{1}$ is 3-periodic with linking number 2. Przytycki and Sokolov \cite{PS} introduced the notion of strongly periodic links as follows.\\

{\bf Definition 2.2} {\sl { Let $p\geq 2$ be an integer, a $p$-periodic link $L$ is said to be strongly periodic if and only if the linking number of each component of $L$ with the axis $\Delta$ is congruent to zero  modulo $p$.}} \\

Now let us  define periodic manifolds:\\
{\bf Definition 2.3} {\sl { Let $p\geq 2$ be an integer,  a 3-manifold $M$ is said to be $p$-periodic if and only if there exists an orientation preserving  auto-diffeomorphism of $M$, such that:\\
1- Fix($h$) is homeomorphic to the  circle $S^{1}$,\\
2- $h$ is of order $p$.\\ }}

The 3-sphere $\S$ is $p$-periodic for all  $p\geq 2$.  The acting diffeomorphism is nothing but the  above defined
 rotation $\varphi_{p}$. Another  concrete example is the Brieskorn manifold $\Sigma(2,3,7)$, which has periods 2, 3 and 7. One can easily explicite the corresponding diffeomorphisms by considering $\Sigma(2,3,7)$ as the intersection of the complex surface $\{(z_{1},z_{2},z_{3} )
 \in \C ^{3}; |z_{1}|^{2}+|z_{2}|^{3}+|z_{3}|^{7}=0\}$ with an appropriate sphere of dimension 5 (see \cite{Mi}).\\

The relationship between periodic links and periodic manifolds
was first noticed by Goldsmith in \cite{Go}. Indeed, she only considered
   cyclic branched covers of $\S$ branched along a knot $K$. In this case, she gave a simple
 algorithm describing how   from the knot $K$ one can construct a  periodic surgery presentation of the covering space. Recently, a much more general result was obtained by Przytycki and Sokolov who showed that any periodic 3-manifold  can be obtained from the 3-sphere by surgery along a strongly  periodic link. To be more precise, they proved the following theorem.\\

 {\bf Theorem 2.4  (\cite{PS}).} {\sl Let $p$ be a prime. A 3-manifold  $M$ is  $p$-periodic if and only if $M$ is obtained from
$\S$ by surgery along a strongly periodic link of period $ p$. }\\

\begin{center}{\sc 3- The quantum $SU(3)$ invariant }\end{center}

{\bf Linear skein theory for $SU(3)$ }\\
This  section is  to review the  construction  of  the quantum $SU(3)$
invariant of 3-manifolds via linear skein theory. This skein theory was
 first introduced by Kuperberg \cite{Ku} in order  to construct
 an invariant of links in the three-sphere. Ohtsuki and Yamada
 \cite{OY} proved that Kuperberg's invariant can be used to define the  $SU(3)$ invariant for 3-manifolds, details about the material exposed in this section can be found in  \cite{Ku}, \cite{MO} and  \cite{OY}.\\
Let $F$ be an oriented surface. The graphs considered in the sequel are  oriented trivalent graphs possibly with loops with no vertices. Moreover, we  assume that at each vertex the three edges are oriented in the same way, all inner or all outer.  A diagram in  $F$ is an oriented trivalent graph immersed in  $F$ such that the singular points are only transverse double points, to each of which over and under crossing information is associated. Equivalently, locally a diagram is a trivalent graph or the diagram of an  oriented link. As usual diagrams are considered here  up to isotopy. Throughout the rest of this paper $\emptyset$ denotes the empty diagram and $\bigcirc$ denotes the trivial knot. The following skein relations were introduced by Kuperberg \cite{Ku}.

\begin{center}
\includegraphics[width=8cm,height=13cm]{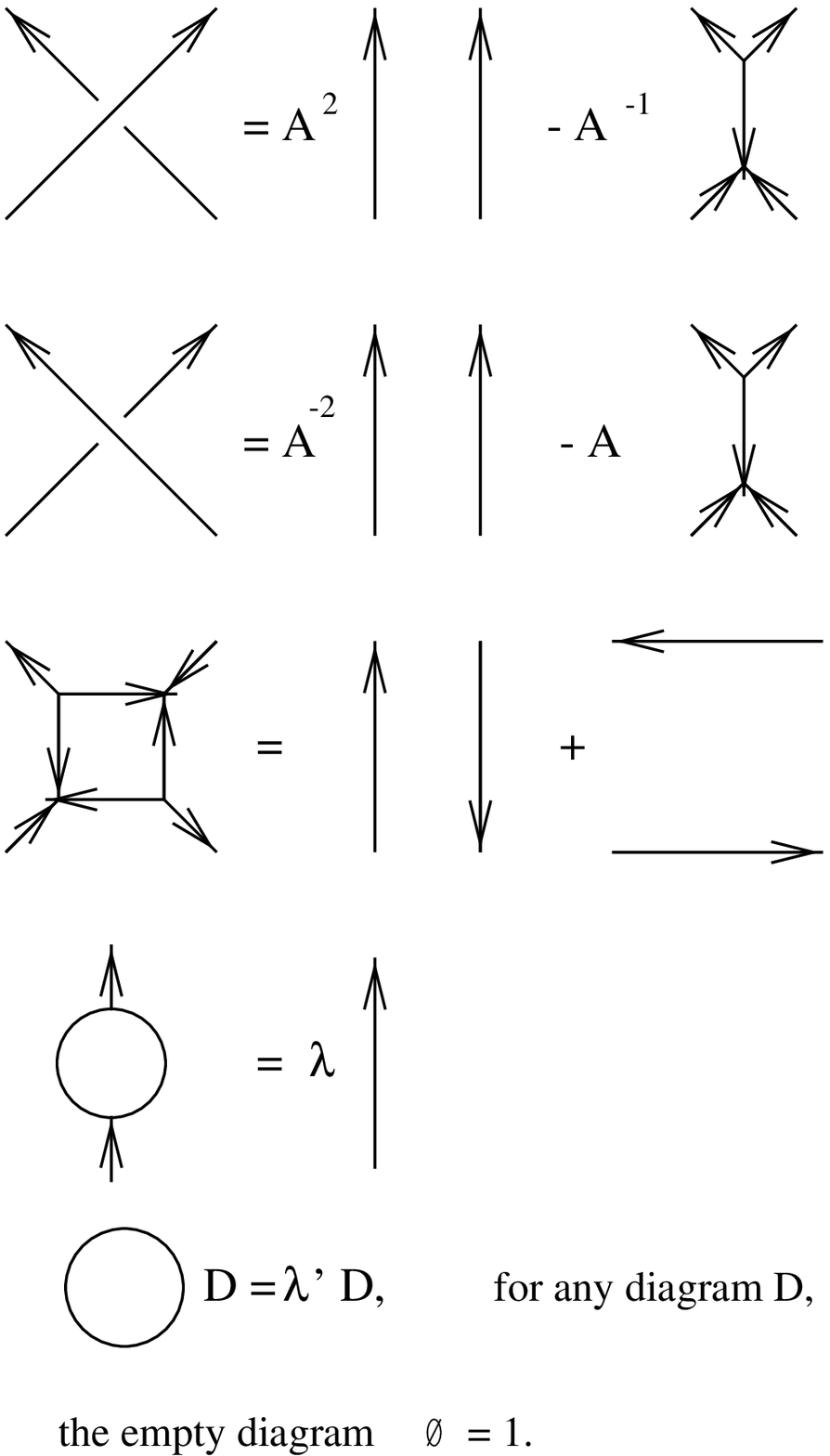}
\end{center}
\begin{center}{\sc  Figure 1. }\end{center}
Where $\lambda =\displaystyle\frac{A^{6}-A^{-6}}{A^{3}-A^{-3}}$ and
$\lambda' =\displaystyle\frac{A^{9}-A^{-9}}{A^{3}-A^{-3}}$.\\

{\bf Remark 3.1} These relations allowed Kuperberg  to define an invariant of oriented  links which we  denote here  by $J$. In fact, this invariant is nothing but a specialization of the Homfly polynomial \cite{MOY}. By setting $q=A^{6}$, this invariant can be defined only by the skein relations:

$$\begin{array}{ll}
{\bf (i)}&J (\emptyset)=1\\
 {\bf (ii)}&J ({\bigcirc}\cup L)=(1+q+q^{-1})J(L), \\
{\bf (iii)}&q^{3/2}J(L_{+})-q^{-3/2}J(L_{-}) =( q^{1/2}-q^{-1/2})J(L_{0})
\end{array}$$
where  $L_{+}$,  $L_{-}$ and  $L_{0}$ are three oriented links which are identical except
 near one crossing where they are as in the following figure:

\begin{center}
\includegraphics[width=8cm,height=2cm]{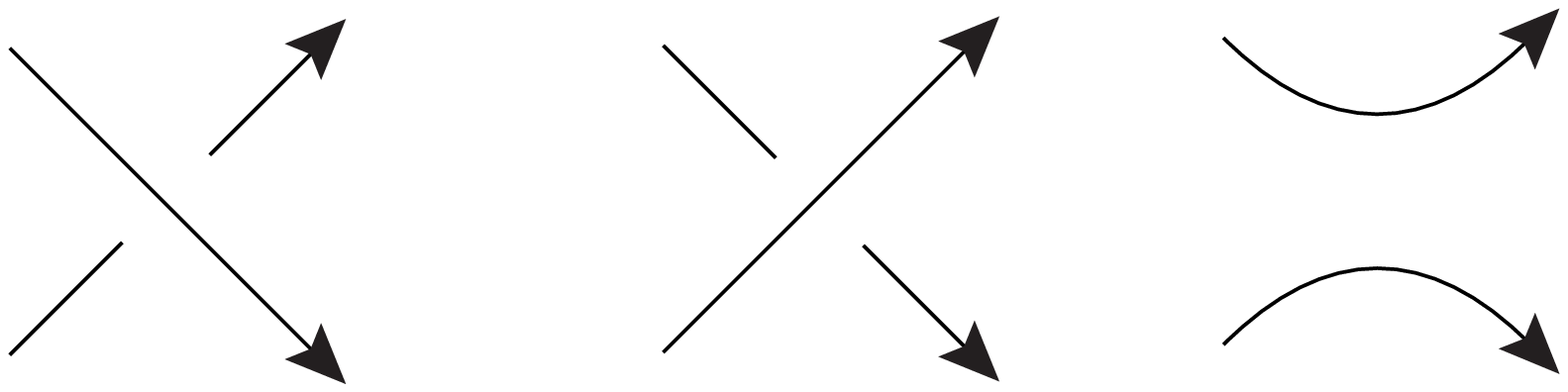}
\end{center}
\begin{center} {\sc Figure 2} \end{center}

{\bf Definition 3.2} {\sl Let $F$ be  an oriented surface. The skein mudule ${\cal S}(F)$ is defined as the quotient
 $\Lambda_{r}$-module of formal linear sums of diagrams in $F$ with coefficient in $\Lambda_{r}$, divided  by the relations of Figure 1.}\\
Let $I$ be  the unit interval $[0,1]$ and  $S^{1}\times I$ the annulus.
There is a natural product on
 ${\cal S}(S^{1}\times I)$. Moreover, if
we consider the two diagrams $x$ and $y$ as in figure 3, we have the following proposition. \\
\null
\vspace{0.5cm}
\begin{center}
\begin{picture}(0,0)
\put(-110,0){$x:$}
\put(-70,0){\circle{10}}
\put(-70,0){\circle{20}}
\put(-80,0){\vector(0,-1){1}}
\put(-70,0){\circle{40}}

\put(30,0){$y:$}
\put(70,0){\circle{10}}
\put(70,0){\circle{20}}
\put(60,0){\vector(0,1){1}}
\put(70,0){\circle{40}}

\end{picture}
\end{center}
\vspace{0.2cm}

\begin{center} {\sc Figure 3} \end{center}

{\bf Proposition 3.3 (\cite{OY}).} {\sl
{\bf (i)} ${\cal S}(\R^{2})$ is isomorphic to  $\Lambda_{r}$.\\
{\bf (ii)} ${\cal S}(S^{1}\times I)$ is isomorphic to $\Lambda_{r}[x,y]$ as an algebra.\\
}
{\bf The invariant ${\cal I}_{r}$.}\\
For any pair of positive integers $(n,m)$, we define the element $P_{n,m}$ as in \cite{OY}. Let $r\geq 4$, we define $\omega _{r}$ as the element of ${\cal S}(S^{1}\times I)$ given by the following formula:
$$\displaystyle\sum_{n+m\leq r-3, n,m\geq 0}[n+1][m+1][n+m+2]P_{n,m}.$$

If $L=l_{1}\cup l_{2}....\cup l_{k}$ is a $k$-component framed link in the three-sphere, we denote by
 $L_{(i_{1},i_{1}'),\dots,(i_{k}, i_{k}')}$ the link obtained from $L$ by decorating the component $l_{j}$ by the element
  $x^{i_{j}}y^{i_{j}'}$ for all $1\leq j \leq k$. We define  the multibracket $\langle .,\dots ,. \rangle_{L}$ as the multilinear form : ${\cal S} \times {\cal S} \times ... \times {\cal S} \longrightarrow\Z[A^{\pm 1}]$ evaluated on set of generators as follows:
$$\langle x^{i_{1}}y^{i_{1}'}, \dots, x^{i_{k}}y^{i_{k}'} \rangle_{L} = J(L_{(i_{1},i_{1}'),\dots,(i_{k}, i_{k}')}).$$

Throughout the rest of this paper if $L$ is a framed link we denote by   $\sigma_{+}(L)$ (resp.  $\sigma_{-}(L)$) the number of positive (resp. negative) eigenvalues of the linking matrix of $L$. In addition, $U_{+}$  (resp $U_{-}$ )  denotes the trivial knot  with framing +1 (resp. -1). Let $M$ be a three-manifold and $L$ a surgery presentation of $M$. The quantum $SU(3)$ invariant is defined by the following formula:
$${\cal I}_{r}(M)= \displaystyle\frac{\langle \omega_{r}, \dots, \omega_{r} \rangle_{L}}{\langle \omega_{r}\rangle_{U_{+}}^{\sigma_{+}(L)} \langle \omega_{r} \rangle_{U_{-}} ^{\sigma_{-}(L)}}.$$
{\bf Remark 3.4}  In \cite{OY}, the authors have noticed that $\langle \omega_{r}\rangle_{U_{+}}$ and
$ \langle \omega_{r} \rangle_{U_{-}}$ are  complex conjugate. Moreover, they proved that $|\langle \omega_{r}\rangle_{U_{+}} |=|\langle \omega_{r} \rangle_{U_{-}}|=\frac{\sqrt{3}r}{|A^{6}-1|^{3}}$. Keeping this in mind, we can see that the only
possible  denominator for  ${\cal I}_{r}$ is a power $3r^{2}$. Consequently,   ${\cal I}_{r}$ may be seen as an element of
$\Lambda_{r}[\displaystyle\frac{1}{3r}]$.

\begin{center}{\sc 4- Proof of theorem 1.1}\end{center}

Several results about the quantum invariants  of periodic links and knots
were obtained in
 \cite{Mu}, \cite{Pr}, \cite{Tr1}, \cite{Tr2}, \cite{Yo1} and \cite{Yo2}. In \cite{Ch2}, we introduced a criterion for periodicity of links using the polynomial invariant $J$. The proof of theorem 1.1 relies heavily on this criterion.  Note that the invariant $J$ corresponds to the invariant $P_{3}$ in \cite{Ch2} and that the variable $q$ corresponds to $A^{6}$.\\
 {\bf Theorem 4.1 (\cite{Ch2})} {\sl Let $p$ be a prime and $L$ a $p$-periodic link. Then we have:
$$J(L) \equiv  (J(\overline  L))^{p} \mbox{ modulo  }p ,[3]^{p}-[3].
 $$}
As in  the case of the Kauffman bracket, the key observation in the proof of theorem 1.1 is that the congruence obtained for the invariant $J$ may be extended to the multibracket of strongly periodic links. Indeed, we can prove the following proposition.\\

{\bf Proposition 4.2} {\sl { Let  $p$ be a prime  and  $L$ a  strongly periodic link with period $p$. Then  for all  $r\geq 4$ the following congruence holds:
$$\langle \omega_{r}, \dots , \omega_{r} \rangle _{L}(A)=(\langle \omega_{r}, \dots , \omega_{r} \rangle _{\overline L}(A))^{p} \mbox{ modulo } p, [3]^{p}-[3].$$}}

Assume proposition 4.2 is proved. Let $M$ be a  $p$-periodic 3-manifold. According to theorem 2.4, $M$ is obtained from the three-sphere by surgery along a strongly periodic link $L$. The quotient $\overline M$ is obtained from $\S$ by surgery along the quotient link $\overline L$. The condition given by proposition 4.2  leads to the following congruence modulo $p$ and   $[3]^{p}-[3]$.
$${\cal I}_{r}(M)\langle \omega_{r}\rangle_{U_{+}} ^{\sigma _{+}(L)} \langle \omega_{r}\rangle_{U_{-}} ^{\sigma _{-}(L)} \equiv
{\cal I}_{r}(\overline M)^{p} \langle \omega_{r}\rangle_{U_{+}} ^{p\sigma _{+}(\overline L)} \langle \omega_{r}\rangle_{U_{+}} ^{p\sigma _{-}(\overline L)}.
$$

From the fact that $|\langle \omega_{r}\rangle_{U_{+}} |=|\langle \omega_{r} \rangle_{U_{-}}|=\frac{\sqrt{3}r}{|A^{6}-1|^{3}}$, we deduce that both
$\langle \omega_{r}\rangle_{U_{+}}$  and $ \langle \omega_{r}\rangle_{U_{-}}$ are invertible in the $\Lambda_{r}[\displaystyle\frac{1}{3r}]$. Thus,

$${\cal I}_{r}(M) \equiv
{\cal I}_{r}(\overline M)^{p} \langle \omega_{r}\rangle_{U_{+}} ^{p\sigma_{+}(\overline L)-\sigma_{+}(L)} \langle\omega_{r}\rangle _{U_{-}}^{p\sigma_{-}(\overline L)-\sigma_{-}(L)}.
$$
As the link $L$ is strongly $p$-periodic, then the number of components of $ L$ is $p$ times the number of components of $\overline L$. Assume that  $\overline L$ has $n$ components. If $M$ is a rational homology sphere then $\overline M$ is also a rational homology sphere and we have the following easy  facts.
$$\begin{array}{rl}
p\sigma_{+}(\overline L)-\sigma_{+}(L)+p\sigma_{-}
(\overline L)-\sigma_{-}(L)&=p(\sigma_{+}(\overline L)
+\sigma_{-}(\overline L))-(\sigma_{+}(L)
+\sigma_{-}(L))\\
&=pn-pn=0.
\end{array}
$$
Thus
$${\cal I}_{r}(M) \equiv{\cal I} _{r}(\overline M)^{p}
(\displaystyle\frac{\langle \omega_{r}\rangle_{U_{+}} }{\langle \omega_{r}\rangle_{U_{-}} })^{p\sigma_{+}(\overline L)-\sigma_{+}
(L)}.
$$
In order to   compute the value of the term $\displaystyle\frac{\langle \omega_{r}\rangle_{U_{+}} }{\langle \omega_{r}\rangle_{U_{-}} }$
we are going to use the calculus developed in \cite{MO} in proving  the invertibility of $\langle \omega_{r}\rangle_{U_{+}}$ and $ \langle \omega_{r}\rangle_{U_{-}}$. In fact, it was established that :
$$\langle \omega_{r}\rangle_{U_{+}}=\displaystyle\frac{-A^{-18}}{(A^{3}-A^{-3})}
(\sum_{r=0}^{r-1}A^{6k^{2}})(\sum_{r=0}^{3r-1}A^{2k^{2}}).$$
Since $\langle \omega_{r}\rangle_{U_{+}}$ and  $\langle \omega_{r}\rangle_{U_{-}}$ are complex conjugate we get :
 $$\langle \omega_{r}\rangle_{U_{-}}=\displaystyle\frac{-A^{18}}{(A^{-3}-A^{3})}
(\sum_{r=0}^{r-1}A^{-6k^{2}})(\sum_{r=0}^{3r-1}A^{-2k^{2}}).$$
Thus,
$$\displaystyle\frac{\langle \omega_{r}\rangle_{U_{+}} }{\langle \omega_{r}\rangle_{U_{-}}}=-A^{-36}
\displaystyle\frac{ (\displaystyle\sum_{r=0}^{r-1}A^{6k^{2}})(\displaystyle\sum_{r=0}^{3r-1}A^{2k^{2}})}{(\displaystyle\sum_{r=0}^{r-1}A^{-6k^{2}})(\displaystyle\sum_{r=0}^{3r-1}A^{-2k^{2}})}.$$
Using the fact that  $|\displaystyle\sum_{r=0}^{r-1}A^{6k^{2}}|^{2}=r$ and $|\displaystyle\sum_{r=0}^{3r-1}A^{2k^{2}}|^{2}=3r$,  we can  conclude that :

$$\displaystyle\frac{\langle \omega_{r}\rangle_{U_{+}} }{\langle \omega_{r}\rangle_{U_{-}}}=-A^{-36}\frac{(\displaystyle\sum_{r=0}^{r-1}A^{6k^{2}})^{2}(\displaystyle\sum_{r=0}^{3r-1}A^{2k^{2}})^{2}}
{3r^{2}}.$$
Setting  $\alpha={p\sigma_{+}(\overline L)-\sigma_{+}(L)}$  ends the proof of theorem 1.1.\\

{\bf Proof of proposition 4.2.}
Here we adapt the techniques  developed in \cite{Ch1} to deal with the multibracket derived from the Kauffman bracket. We shall  briefly recall the main idea of the proof and  refer the reader to \cite{Ch1} for details. Let $L$ be  a framed link, we denote by
${\cal E}_{L}$ the set of all possible  decorations of $L$ coming from $\omega _{r}$. The element
$\omega _{r}$ is of the form  $\displaystyle \sum_{i \geq 0} Q_{i,i'}(A)x^{i}y^{i'}$ for some polynomials
 $Q_{i,i'}(A)$. Therefore, the multibracket is given by  the following formula:

$$\langle \omega_{r}, \dots , \omega_{r} \rangle _{L}=
\displaystyle \sum_{K\in {\cal E}_{L}} Q_{K}(A) J(K).
$$
Let $L$ be  a stongly $p-$periodic link and $\overline L$ its factor link. Assume that  $\overline L=l_{1}\cup l_{2} \cup ... \cup l_{n}$. As the linking number of each component $l_{i}$ with the axis $\Delta$ is zero modulo $p$ then the link $L$ is of the following  form:
$$(l_{1}^{1}\cup
\dots \cup l_{1}^{p})\cup (l_{2}^{1}\cup \dots \cup l_{2}^{p})\cup
\dots\cup (l_{n}^{1}\cup \dots \cup l_{n}^{p}).$$
Moreover, for all $i$ the components $l_{i}^{j}$ are identical and   cyclically permuted  by the rotation $\varphi _{p}$.

{\bf Definition 4.3} {\sl  Let $L$ be a strongly periodic link with period $p$. A decoration of $L$ is said to be  periodic  if and only if for all $i$  the components  $l_{i}^{j}$ are decorated in the same way. Let   ${\cal E}'_{L}$ be the subset of ${\cal E}_{L}$ made up of periodic decorations}.\\
{\bf Lemma 4.4} {\sl
Let $p$ be a prime and $L$ a strongly $p$-periodic link then we have the following:\\
$$\langle \omega_{r}, \dots , \omega_{r} \rangle _{L}\equiv
\displaystyle \sum_{K\in {\cal E'}_{L}} Q_{K}(A) J(K) \mbox{ modulo } p.
$$}

Proof of lemma 4.4. The proof  is based on the fact that the cyclic group $\Z/p\Z$ acts on ${\cal E}_{L}$. As $p$ is prime, orbits of this  action are made up of 1 or $p$ elements. In the case of a  $p$-elements orbit, these elements are identical and their contribution to the multibracket vanishes  modulo $p$. A 1-element orbit corresponds to a periodic decoration of $L$. Hence, for the computation of the multibracket modulo $p$ one may only consider elements of ${\cal E}'_{L}$, see \cite{Ch1} for details.\\
Now note that if $K$ is in  ${\cal E}'_{L}$ then $K$ is of the form :
$$K=L_{(i_{1},i_{1}'),\dots,(i_{1},i_{1}'),(i_{2},i_{2}'),\dots,(i_{2},i_{2}'), \dots ,(i_{k}, i_{k}'),\dots,(i_{k}, i_{k}')}
$$
where each  $(i_{j}, i_{j}')$ appears $p$-times. Thus $K$ is a $p$-periodic link and its quotient link is:
$$\overline K=\overline L_{(i_{1},i_{1}'),(i_{2},i_{2}'), \dots,(i_{k}, i_{k}')}.$$
Therefore for all $K \in {\cal E}'_{L}$, the condition given by  of theorem 4.1 is satisfied :
$$J(K) \equiv  (J(\overline K))^{p} \mbox{ modulo  }p ,[3]^{p}-[3].
 $$
Finally, it is easy to see that  every element of
${\cal E}'_{L}$ induces a natural decoration of ${\overline L}$, thus an element of  ${\cal E}_{\overline L}$, and vice-versa. Moreover the coefficients that appear in the sum satisfy  $Q_{K}(A)=(Q_{\overline K}(A))^{p}$. This ends the proof of proposition 4.2.\\

\begin{center}{\sc 5- Similar results for the  MOO invariant}\end{center}

In \cite{MOO}, Murakami, Ohtsuki and Okada introduced an invariant of 3-manifolds so called the MOO invariant. This invariant is defined from the linking matrix of the surgery presentation of the 3-manifold. Let $M$ be a 3-manifold obtained from  $\S$ by surgery  along an $m$-component  framed link $L$. By ${\cal B}_{L}$ we denote the linking matrix of $L$. Let $A$ be a primitive root of unity of order $N$ (resp. $2N$) if $N$ is odd (resp. even). Let  $G_{N}(A)$ denotes the Gaussian sum
$\displaystyle\sum_{k \in \Z/N\Z}A^{k^{2}}$. The MOO invariant is defined as follows:

$$
 Z_{N}(M)=(\displaystyle\frac
{G_{N}(A)}{|G_{N}(A)|})^{\sigma(L)}|G_{N}(A)|^{-m}\displaystyle\sum_{l \in (\Z/N\Z)^{m}}A^{t(l)
{\cal B}_{L}l}.
$$
Where $\sigma(L)$ is the signature of the linking matrix ${\cal B}_{L}$, $l$ is regarded as a column vector and $t(l)$ is its transposed row vector.
As in the case of the other quantum invariants, the MOO invariant can be constructed via skein theory \cite{Gil}. We shall  recall here briefly this
very simple construction. Let ${\Lambda}'_{N}=\Z[A^{\pm 1}]_{/\Phi_{N}(A)}$. We denote by ${\cal L}(M)$ the free  ${\Lambda}'_{N}$-module generated by the set of links in the 3-manifold $M$. Now consider the following skein relations:
$$\begin{array}{rl}
{\cal R}_{1}:& \emptyset=1\\
{\cal R}_{2}:& L_{+}=A L_{0}\\
{\cal R}_{3}:&L_{-}=A^{-1} L_{0}\\
 {\cal R}_{4}:& L\cup \bigcirc= L

\end{array}
 $$
The linking skein module of $M$ which is  denoted here as  ${\cal S}'(M)$  is defined as the quotient module of ${\cal L}(M)$ by relations ${\cal R}_{i}$ for  $1\leq i \leq 4$. We can easily prove that
${\cal S}(\S)={\Lambda}'_{N}$ and that  the linking module of the solid  torus   ${\cal S}'(S^{1}\times D^{2})$ is isomorphic as an algebra to
${\Lambda}'_{N}[x,x^{-1}]$ where $x$ and $x^{-1}$ are as in figure 4.\\
\null
\vspace{0.5cm}
\begin{center}
\begin{picture}(0,0)
\put(-110,0){$x:$}
\put(-70,0){\circle{10}}
\put(-70,0){\circle{20}}
\put(-80,0){\vector(0,-1){1}}
\put(-70,0){\circle{40}}

\put(20,0){$x^{-1}:$}
\put(70,0){\circle{10}}
\put(70,0){\circle{20}}
\put(60,0){\vector(0,1){1}}
\put(70,0){\circle{40}}

\end{picture}
\end{center}
\vspace{0.2cm}

\begin{center} {\sc Figure 4} \end{center}

If $L$ is an $m$-component  oriented  link, we define the bracket of $L$ by $\langle L \rangle = A^{lk(L)}$, where $lk(L)$ is the total linking number of $L$. The multibracket  $\langle .,\dots, . \rangle_{L}$ is  the multilinear form ${\cal S}' \times {\cal S}' \times ... \times {\cal S}' \longrightarrow \Lambda'_{r}$   defined on the set of generators by
  $\langle x^{i_{1}},x^{i_{2}},\dots ,x^{i_{m}} \rangle_{L}=\langle L_{i_{1},i_{2},\dots ,i_{m}} \rangle$. Let $N \geq 1$ be  an integer and $\omega_{N}=\displaystyle\sum_{k=0}^{N-1}x^{k}$. The MOO  invariant is obtained as follows:
$$Z_{N}(M)=|G_{N}(A)|^{-b_{1}(M)} \displaystyle\frac{\langle\omega_{N}, \dots ,\omega_{N} \rangle_{L}}{ \langle \omega_{N}\rangle_{U_{+}}^{\sigma_{+}(L)} \langle \omega_{N} \rangle_{U_{-}} ^{\sigma_{-}(L)}}.$$
Where $b_{1}(M)$ is the first  Betti number of the manifold $M$. In the following we only consider the case $N$ is odd. We are going to consider the MOO invariant as an element of $\Lambda'_{r}[\displaystyle\frac{1}{N}]$ (this is an immediate consequence from the defintion).  Let $M$   be a $p$-periodic rational homology sphere. We know that $M$ can be  obtained from the three-sphere by surgery along a strongly $p$-periodic link  $L$. Using the easy fact that $\langle L \rangle=\langle \overline L \rangle^{p} $ and the techniques of paragraph 4 we can prove that :
$$\langle \omega_{N}, \dots , \omega_{N} \rangle _{L}(A)=(\langle \omega_{N}, \dots , \omega_{N} \rangle _{\overline L}(A))^{p} \mbox{ modulo } p.$$
This means that for $N$ coprime to $p$ we have :
$$
Z_{N}(M)=(\displaystyle\frac{\langle \omega_{N}\rangle_{U_{+}}}{\langle \omega_{N}\rangle_{U_{-}}})^{p\sigma_{+}(\overline L)-\sigma_{-}(L)}
(Z_{N}(M))^{p} \mbox { modulo p}.
$$
It is easy to see that $ \langle \omega_{N}\rangle_{U_{+}}=\overline{\langle \omega_{N}\rangle_{U_{-}}}=G_{N}(A).$ As $N$ is odd, elementary properties of Gaussian sums \cite{La} shows that:\\
 $$\displaystyle\frac{\langle \omega_{N}\rangle_{U_{+}}}{\langle \omega_{N}\rangle_{U_{-}}}=\displaystyle\frac{G_{N}(A)}{\overline{G_{N}(A)} }=\pm 1.
$$
Thus we get the following:\\
{\bf Theorem 5.1} {\sl Let $p$ be  a prime and $M$ a rational homology 3-sphere. If $M$ is $p$-periodic then for all odd
 $N\geq 2$ coprime to $p$, the following congruence holds:
$$ {Z }_{N}(M)\equiv \pm ({Z}_{N}(\overline M ))^{p} , \mbox{ modulo }p $$
where $\alpha$ is an integer.}\\

\begin{center}{\sc 6- Concluding remarks}\end{center}
{\bf Remark 6.1} The key observation in the proof of our main result here is the condition given by theorem 4.1. Recently, Przytycki
and Sikora \cite{PSi} obtained a generalisation of this condition for  the $SU(N)$ invariants of periodic links, $N\geq$
is an odd integer. On the other hand, Yokota \cite{Yo} introduced a skein theory for the $SU(N)$ quantum  invariants of 3-manifolds. This
indicates that the methods used in our paper may lead to a generalisation of theorem 1.1 to the  $SU(N)$ quantum  invariants of periodic rational homology spheres. We plan to describe this in a forthcoming paper.\\
{\bf Remark 6.2} On the lines of the ideas discussed in the present paper. We discuss in \cite{Ch3}, the surgery formula introduced by  Lescop
\cite{Le} for the Casson-Walker invariant. Namely, we prove a necessary condition for the Casson invariant of periodic three-manifolds.


\begin{footnotesize}
{\sc Permanent Address:}\\
Nafaa Chbili\\
 D\'epartement de Math\'ematiques, Facult\'e des Sciences de Monastir.\\
Boulevard de l'environnement,  \\
Monastir 5000, Tunisia.\\
e-mail: {\underline {nafaa.chbili@esstt.rnu.tn}}

\end{footnotesize}

\end{document}